\documentclass[twoside]{article}

\usepackage{amsmath,amsthm,amsfonts,latexsym,amscd,amssymb,enumerate}
\usepackage[hypertex]{hyperref}

\swapnumbers
\theoremstyle{plain}
\newtheorem{theorem}{Theorem}[section]
\newtheorem{lemma}[theorem]{Lemma}
  
\newtheorem{proposition}[theorem]{Proposition}

\newtheorem{assumption}[theorem]{Assumption}
\theoremstyle{definition}

\newtheorem{set-up}[theorem]{Geometric set-up}

\DeclareMathOperator{\Ind}{Ind}

\newcommand{\forget}[1]{}

\def  \nuint {\raise10pt\hbox{$\nu$}\kern-6pt\int}

\newcommand\ep{\operatorname{\epsilon}}

\def \P{\mathcal P}

\def \Sp {{\cal S}}

\newcommand\ha{\frac12}

\newcommand\Di{D\kern-6pt/}
\newcommand\cDi{{\mathcal D}\kern-6pt/}
\newcommand\spi{S\kern-6pt/}
\newcommand \cspi{\Sp\kern-6pt/}

\newcommand\CC{\mathbb C}

\def \cal {\mathcal}

\newcommand\AAA{\mathbb A}

\newcommand\FF{\mathbb F}

\newcommand\pa{\partial}

{\catcode`@=11\global\let\c@equation=\c@theorem}

\allowdisplaybreaks[2]

\title{The index of Dirac operators on manifolds with fibered boundaries}

\author{Eric Leichtnam \\ Institut de Math\'ematiques de Jussieu\\175 rue de Chevaleret\\75013 Paris, France\and
Rafe Mazzeo \\ Department of Mathematics\\
Stanford University\\
Stanford, CA 94305 USA
\and and Paolo Piazza \\Dipartimento di Matematica\\
Universit\`a di Roma  "La Sapienza"\\P.le A. Moro 2\\00185 Roma, Italy} 

\begin{document}

\maketitle

\begin{abstract}
Let $X$ be a compact manifold with boundary $\pa X$, and suppose
that $\pa X$ is the total space of a fibration 
\[
Z\rightarrow \pa X \rightarrow Y\, .
\]

Let $D_\Phi$ be a generalized Dirac operator associated
to a $\Phi$-metric $g_\Phi$ on $X$. Under the assumption that
$D_\Phi$ is fully elliptic we prove an index formula for $D_\Phi$.
The proof is in two steps: first, using results of Melrose and Rochon, 
we show that the index is unchanged if we  pass to a certain
$b$-metric $g_b (\epsilon)$. Next we write the $b-$ (i.e.\ the APS) index formula
for $g_b(\ep)$; the $\Phi$-index formula follows by analyzing the limiting behaviour 
as $\epsilon\searrow 0$ of the two terms in the formula. The interior term
is studied directly whereas the adiabatic limit formula for the eta invariant
follows from work of Bismut and Cheeger.
\end{abstract}



\thanks{ 

eta invariant}
\section{Introduction}\label{sec:intro}

Let $X$ be an even dimensional, compact, oriented spin manifold with boundary such that 
$\pa X$ is the total space of a fibration $Z^\ell \rightarrow \pa X \xrightarrow{\phi} Y^k.$ 
(Thus $\dim X=\ell+k+1=2m$.) There are many interesting index formul\ae\ for twisted
Dirac operators $D$ on $X$ corresponding to various different classes of complete 
metrics $g$ on the interior of $X$. Under certain hypotheses which ensure that $D$ 
is either Fredholm, or at least has finite $L^2$ index, and that the usual Atiyah-Singer
density has finite integral, the goal is to identify the index defect, i.e.\ the 
difference between $\Ind (D)$ and the Atiyah-Singer integrated characteristic form. Most prominent,
of course, is the Atiyah-Patodi-Singer theorem when $g$ has asymptotically cylindrical ends,
in which case the index defect is (minus one half) the eta invariant of the induced twisted
Dirac operator on $\pa X$ \cite{Melrose}. This does not take advantage of the fibred boundary structure.
Two interesting classes of metrics which do take this into account are the
fibred boundary and fibred cusp metrics, also called $\Phi-$ and $d-$ metrics,
respectively. These appear naturally in many interesting geometric settings, cf.\ \cite{HHM}:
for example, complete Ricci flat metrics are often $\Phi-$metrics, while
locally symmetric metrics with ${\mathbb Q}-$rank one cusps are $d-$metrics.

To define these, introduce the following notation. Fix a splitting $T(\pa X) = 
T_V(\pa X)\oplus T_H (\pa X)$ into vertical and horizontal subspaces, where 
$T_V(\pa X) = T(\pa X/Y)$ is the fibre tangent bundle, and $T_H (\pa X)$ is 
identified with $\phi^*(TY)$. We consider metrics $\tilde{g}$ 
on $\pa X$ and $h$ on $Y$ so that $\phi$ is a Riemannian submersion. This means that the 
restriction of $\tilde{g}$ to $T_H(\pa X)$ is identified with $\phi^* h$, and the 
subbundles $T_H(\pa X)$ and $T(\pa X/Y)$ are orthogonal. We write $\tilde{g} = \phi^*h + 
\kappa$, where $\kappa$ is a symmetric two-tensor on $\pa X$ which is positive definite 
on $T(\pa X/Y)$. 

Let $x$ be a defining function for $\pa X$ in some neighbourhood of the boundary. Suppose
also that $h$ and $\kappa$ are allowed to depend smoothly on $x$, all the way to $x=0$. Then
an exact $b$-metric and an exact cusp ($c$-) metric on $X$ are ones which have the form
\[
\frac{dx^2}{x^2} + \tilde{g}, \qquad \frac{dx^2}{x^4} + \tilde{g}
\]
in this neighbourhood, respectively ; likewise, exact $\Phi-$ metrics and exact $d-$metrics have the forms
\[
\frac{dx^2}{x^4} + \frac{\phi^* h}{x^2} + \kappa, \qquad \mbox{and} \qquad
\frac{dx^2}{x^2} + \phi^* h + x^2 \kappa,
\]
respectively, in this neighbourhood. 
(The term `exact' in each of these refers to the fact that there
are no cross-terms, at least to principal order; this is a natural, but not a serious assumption, 
and there are generalizations of the ideas and formul\ae\ we discuss here to the various
`nonexact' settings.) For simplicity in all of the discussion below, we usually label
a metric as $g_b$, $g_c$, $g_\Phi$ and $g_d$ to indicate that it is one of these four types.
Note also that when discussing $b-$ and $c-$metrics, it is not  important that
$\tilde{g}$ respect the fibration structure (nor, of course, even that $\pa X$ have such
a structure). 

Assume that $X$ and $Y$ are spin, and fix spin structures on each of these manifolds; there is
an induced spin structure on the fibres $\phi^{-1}(Y) := Z_y \subset \pa X$. The ($\mathbb{Z}_2-$graded)
spin bundles on $X$, $\pa X$ and $Z_y$ are denoted $S$, $S^{\pa}$ and $S^{Z_y}$, respectively. 
Let $E\rightarrow X$ be an hermitian complex vector bundle endowed with a unitary connection. 
Fixing also a metric $g$ of any of the types above, we obtain a twisted Dirac operator 
\[
D_g^+: C^\infty(X, E\otimes S^+) \rightarrow C^\infty(X, E\otimes S^-).
\]
If $g$ is of one of the preceding types, then we also write $D_b$, $D_c$, $D_\Phi$, $D_d$
for the corresponding Dirac operator to indicate its asymptotic type.
The associated boundary operator $D^{\pa}$ induces a family of Dirac operators 
$\{D^\pa_y\}_{y\in Y}$, where each $D^\pa_y$ acts on $C^\infty(\phi^{-1}(y), E \otimes S^\partial_y)$.

For $b-$ and $c-$ metrics, the simplest form of the APS occurs when $D^\pa$ is invertible,
although the general result if this is not satisfied is not much more difficult. In the
other two settings, however, the analogous hypothesis is the
\begin{assumption}
\label{assumption}  
For some $\delta>0$, 
\begin{equation}
\label{eq:assumption}
{\rm spec} (D^{\pa}_y) \cap (-\delta,\delta)=\emptyset,  \quad  \forall\,y\in Y
\end{equation}
\end{assumption}
\medskip

The index formula for $D_g$ is known when $g$ is a metric of type $b$, $c$ or $d$; as 
already noted, the first of these is just the APS theorem, while the second in fact
reduces to this theorem in a rather simple way. (This is proved below.) 
The index formula for $d-$metrics is due in the special case of locally symmetric metrics
to M\"uller \cite{Mueller}, and in this general geometric setting was accomplished 
by Vaillant \cite{Vaillant-phd}.  The index defect in this case is the integral 
over $Y$ of the Bismut-Cheeger eta form. (Actually, Vaillant's result holds under
the weaker hypothesis that $\ker D^\pa_y$ has constant rank, in which case the
index formula has an additional boundary contribution.)

Assuming (\ref{assumption}), $D_\Phi^+$ is a fully elliptic operator in the (pseudo)differential $\Phi-$calculus 
developed in \cite{Ma-Mel-phi} and \cite{Vaillant-phd}, and the parametrix construction 
there shows that $D_\Phi^+$ is Fredholm acting between the appropriate ($\Phi$-) Sobolev spaces. 
Answering a question raised in \cite{Ma-Mel-phi}, we prove here that
\begin{theorem} Assuming (\ref{assumption}), and using the notation above, we have
\begin{equation}\label{main-formula} 
\Ind (D_\Phi^+)=\int_X \widehat{A}(X,g_\Phi) \wedge {\rm Ch}\, E  -\ha\int_Y  
\widehat{A}(Y,h)\wedge \widetilde{\eta},
 \end{equation}
where $\widetilde{\eta} \in \Omega^*(Y)$ is the Bismut-Cheeger eta form \cite{BC} for the 
boundary family $(D^{\pa}_y)_{y\in Y}$. 
\end{theorem}

While it is likely that the index formula for this operator can be obtained by
methods similar to those employed in  \cite{Vaillant-phd} for $d$-metrics, that 
proof is very long and
difficult, and it is a reasonable goal to obtain this formula as a consequence
either of that theorem or of the APS theorem.

The equality of the $\Phi$-index and the $d$-index (when  the boundary family
          is invertible) has been  recently proved by Sergiu  Moroianu  \cite{phi-vs-di}
          by reducing
          it directly to Vaillant's theorem \cite{Vaillant-phd}: if $g_d = x^2 g_ \Phi$, and both are exact,
          then $\mbox{Ind}\,(D_\phi) = \mbox{Ind}\, (D_d)$. By Vaillant \cite{Vaillant-phd},
          \[
\Ind (D_d) = \int_X \widehat{A}(X,g_d) \wedge {\rm Ch}\, E  -\ha\int_Y 
\widehat{A}(Y,h)\wedge \widetilde{\eta},
\]
          so it suffices to show that the first integral on the  right is the same as the
          corresponding one for $g_\Phi$, i.e. that
         \begin{lemma}\label{d-integral=phi-integral}
$\int_X \widehat{A}(X,g_\Phi)  \wedge {\rm Ch}\, E =\int_X \widehat{A}(X,g_d)  \wedge {\rm Ch}\, E.$
\end{lemma}

%

Notice  that the two integrals are well defined: this is discussed in 
 \cite[Section 1]{Vaillant-phd}.
Lemma \ref{d-integral=phi-integral} follows simply because $\widehat{A}(X,g_\Phi) = \widehat{A}(X,g_d)$ pointwise,
by conformal invariance. Note too that by a standard transgression argument,
the integrals are equal even when $g_d$ and $x^2 g_\Phi$ coincide only in
a neighbourhood of $\pa X$. 

The proof of (\ref{main-formula}) here is indirect too, but it involves only a reduction 
to the much simpler APS theorem. We shall use the technique
of adiabatic limit, as described below. We first deform $g_\Phi$ to a $b$-metric 
$g_b (\ep)$. The index is unchanged through this deformation, and hence equals the 
index of the Dirac operator corresponding to $g_b(\ep)$. This follows from the analysis 
of Melrose and Rochon \cite{mel-rochon}, specifically their construction of parametrices 
which are uniform in an adiabatic parameter $\ep$ for certain parts of this metric deformation. 
In the (APS) index formula for this $b-$ metric we then take the limit as $\ep \to 0$. 
The fact that the eta invariant term has the correct limiting behaviour follows from the 
Bismut-Cheeger theory \cite{BC} so it remains only to analyze the limiting behaviour 
of the interior integral, which is the new calculation here. The particular metric
family $g_b(\ep)$ is chosen because the Atiyah-Singer integrand for it has the best
behaviour in the limit.

In the initial stages of our work, the plan was to develop a more direct deformation
connecting $g_\Phi$  and $g_b$ and to use a parametrix method to analyze this adiabatic 
limit. However, just at this time the paper of Melrose and Rochon \cite{mel-rochon} appeared, and 
Lemma C.1 there (i.e. Lemma \ref{lem:mel-rochon}  below) allowed us to develop the 
particular and much shorter route presented here. By relying on their substantial and 
deep work, as well as that of \cite{BC}, we are able to give a fairly quick proof of this index formula.

It should be possible, and  would still be of genuine interest, to prove the $\Phi-$index 
theorem directly using heat equation methods. In particular, one would hope to obtain 
another derivation of the fundamental Bismut-Cheeger result in the course of this.

We conclude this discussion by noting that Lauter and Moroianu \cite{LM} prove formula 
\eqref{main-formula} in the special case $Y=S^1$. In fact, in their earlier
paper \cite{LM0}, they also treat the case where $Y$ is arbitrary and establish a less 
precise index formula using homological methods based on ideas of Melrose-Nistor. 
We refer also to \cite{Nye-Singer} for a related formula when $\phi: \partial 
X =S^1 \times S^2 \rightarrow S^2.$
 
We shall prove formula \eqref{main-formula} assuming that $E$ is the trivial line bundle 
$ X \times \CC.$  This is for notational simplicity only, and the general formula
may be deduced using exactly the same reasoning.  In the next section we introduce
the sequence of metric homotopies and prove that the index is unchanged under these
deformations. In the third section we analyze the other side of the index formula,
and especially its behaviour in the adiabatic limit.

\medskip
\noindent
{\bf Acknowledgments.}
We wish to thank Richard Melrose and Fr\'ed\'eric Rochon for explaining their work to us,
and also Sergiu Moroianu for making some valuable suggestions on an early draft of
this note. This work was initiated during a visit by the third author to
Stanford University, and he wishes to thank that department for its hospitality. The research of
Eric Leichtnam
and Paolo Piazza is  partially supported  by a CNR-CNRS  bilateral project.
Rafe Mazzeo was supported by the NSF grant DMS-0505709.

\section{Reduction of $\Ind(D_\Phi)$ to $\Ind(D_b)$}\label{sec:phi-to-b}

In order to avail ourselves of the work of Melrose and Rochon, the homotopy of
metrics we consider consists of the following steps: first deform $g_\Phi$ to
the cusp metric
\begin{equation}\label{eq:c-metrics-1}
g^1_c (\ep):=
\frac{dx^2}{x^4} + \frac{\phi^* h}{(x+\ep)^2}+\kappa;
\end{equation}
next, deform $g^1_c(\ep)$ to the cusp metric
\begin{equation}\label{eq:c-metrics-bis}
g^0_c (\ep):=
\frac{dx^2}{x^4} + \frac{\phi^* h}{\ep^2}+\kappa;
\end{equation}
from here deform in succession to the following three $b$-metrics:
\begin{eqnarray}
g_b^0 (\ep) & := & \frac{dx^2}{x^2} + \frac{\phi^* h}{\ep^2}+\kappa \label{eq:b-metrics-bis} \\
g^1_b (\ep) & := & \frac{dx^2}{x^2} + \frac{\phi^* h}{(x+\ep)^2}+\kappa \label{eq:b-metrics-1} \\
g^2_b(\ep) &:=& { (d x)^2 \over x^2 (x + \epsilon)^2} + { \phi^* h \over (x+\epsilon)^2} + \kappa \label{eq:b-metrics-t}
\end{eqnarray}

Of course we have only specified the forms of these metrics in a fixed collar neighbourhood
of $\pa X$, but we can extend these to the interior arbitrarily, and standard results show
that neither their indices nor the integrals depend on these extensions. 

We denote by $D_*^j(\epsilon)$, $*=c,b$ and $j=0,1,2$, the Dirac operators associated to these 
metrics, respectively.

The first main fact is the
\begin{lemma}\label{lem:b-invertibility}
Assuming (\ref{assumption}), then each of the operators $D_*^j(\ep)$ is fully elliptic
when $\ep>0$ is sufficiently small. 
\end{lemma}
Full ellipticity in either the $b-$ or $c-$ pseudodifferential calculi is simply 
the assumption that not only the interior symbol but also the boundary `indicial operator' 
is invertible.  This follows from Theorem (4.41) in \cite{BC} when $\ep$ is small. Using 
the full ellipticity, one may construct parametrices modulo compact remainders in the appropriate
pseudodifferential calculi. Hence each of the operators $D^j_*(\ep)$ is
Fredholm on the appropriate geometric Sobolev spaces. 

From now on we shall omit mention that the hypothesis (\ref{assumption}) is always
in force here. Furthermore, we shall always assume that $0 < \ep < \ep_0$ for
some sufficiently small $\ep_0$.

We deform to $g_b^2(\ep)$, rather than any of the `simpler' $b$-metrics because
this is the metric for which we can more effectively analyze the limit of
the Atiyah-Singer integrand as $\ep \searrow 0$.

We now present a series of lemmata which state that the indices of the Dirac operators
remains the same through this entire deformation.

The first step uses the work Melrose and Rochon and is the most serious one analytically.
\cite{mel-rochon}. 

\begin{lemma}\label{lem:mel-rochon}
$\Ind(D_\Phi) = \Ind(D^1_c (\epsilon))$.
\end{lemma}

\begin{proof} In Appendix C of \cite{mel-rochon}, Melrose and Rochon consider 
an adiabatic metric deformation connecting a $\Phi$ metric to a $c$-metric. 
Actually, they consider a slightly more general situation where $\pa X$ is the 
total space of a tower of fibrations $\pa X \to 
\tilde{Y} \to Y$, and a corresponding transition between a $\Phi-$metric associated
to the first fibration and a $\Phi-$metric associated to the second. 
By constructing parametrices in an adiabatic calculus, they prove in 
Proposition C.1 of \cite{mel-rochon} that the indices of the Dirac operators
associated to the metrics in this family remain invariant in this passage 
to an adiabatic limit. This assumes that the `boundary symbols,\ i.e.\ 
the normal operators $\mbox{ad}(P)$ and $N(P)$ are invertible, which 
follows directly from our hypothesis (\ref{assumption}). 
\end{proof}

\begin{lemma}\label{lem:c-homotopy} We have
$\Ind(D^1_c (\epsilon)) = \Ind(D_c^0(\epsilon))$, 
$\Ind(D_b^0(\epsilon)) = \Ind(D_b^1 (\epsilon))$ 
and $\Ind(D_b^1 (\ep)) = \Ind(D_b^2(\ep))$.
\end{lemma}

\begin{proof} In each case we simply follow the obvious homotopy of metrics. Thus,
for the cusp setting, let
\begin{equation}\label{eq:homotopy}
g_c(t,\ep):= \frac{dx^2}{x^4} + \frac{\phi^* h}{(tx+\ep)^2}+\kappa\,,\qquad 0 \leq t \leq 1,
\end{equation}
so that $g_c(0,\ep)= g_c^0 (\ep)$, $g_c(1,\ep)= g_c ^1(\ep)$.  The indicial family 
of the corresponding Dirac operators $D_c(t,\ep)$ is independent of $t$, and hence 
each $D_c(t,\ep)$ is Fredholm, so the index is constant.  The argument in the other
two cases is the same. 
\end{proof}

\begin{lemma}\label{lem:homotopy}
$\Ind(D_c^0 (\epsilon)) = \Ind(D_b^0 (\epsilon))$.
\end{lemma}
\begin{proof} As in \cite{mel-rochon}, Lemma (14.1), $\Ind(D_c^0(\ep)$ equals
the index for the incomplete metric $du^2 + \phi^* h/\ep^2 + \kappa$ with APS 
boundary conditions. Since the boundary operator is invertible, this also
equals $\Ind(D_b^0 (\epsilon))$.
\end{proof}

Taken together, this chain of equality gives the
\begin{proposition}\label{propo:phi=b}
$\Ind (D_\Phi)=\Ind (D_b^2(\epsilon))$. 
\end{proposition}

\section{The adiabatic limit}\label{sec:adiabatic}

At this point we simply notation and simply write $g(\ep)$ instead of $g_b^2(\ep)$ .

We begin with the
\begin{proposition}
Assuming, as always, that (\ref{assumption}) holds, then for $\ep$ sufficiently small, 
\begin{equation}
\Ind(D_\Phi)=\int_X \mbox{AS}(g(\ep)) -\frac{1}{2}\eta(D^\pa_{g(\ep)}). 
\label{eq:phigb2}
\end{equation}
\end{proposition}

\begin{proof} Define $\xi = x/(x+\ep)$, so that $d\xi/\xi = \ep dx/x(x+\ep)$. In terms
of this new boundary defining function, $g(\ep) := \ep^{-2} \hat{g}$, where 
\[
\hat{g} = \frac{d\xi^2}{\xi^2} + \ep^2\left(\frac{\phi^*h}{(x+\ep)^2} + \kappa\right). 
\]
The middle term on the right has been kept expressed in terms of $x$ simply to emphasize that 
$\hat{g}$ is an exact $b$-metric which induces $\ep^2$ times the boundary metric induced 
by $g(\ep)$. 

Applying the usual APS formula to $\hat{g}$ gives
\begin{equation}
\Ind(D_{\hat{g}}) = \int_X \mbox{AS}(\hat{g}) - \frac{1}{2}\eta(D^\pa_{\hat{g}}).
\label{eq:apshatg}
\end{equation}
However, clearly  $D_{g(\ep)}$ has the same index as $D_{\hat{g}}$, which is then
the same as $\Ind(D_\Phi)$. Furthermore, using the fact that $\hat{g}$ and 
$g(\ep)$ differ by a constant, we get both
\[
\int_X \mbox{AS}(\hat{g}) = \int_X \mbox{AS}(g(\ep)), \qquad \mbox{and}\qquad \eta(D^\pa_{\hat{g}}) = 
\eta(D^\pa_{g(\ep)}).
\]
Replacing each term in (\ref{eq:apshatg}) with the corresponding quantity 
for $g(\ep)$ gives (\ref{eq:phigb2}).
\end{proof}

The final steps of the proof of the main theorem consist in analyzing the
limiting behaviour as $\ep \searrow 0$ of the two terms on the right in (\ref{eq:phigb2}).

\subsection{Limiting behaviour of the integrand}

\begin{proposition}\label{prop:limit-of-as}
The integral of the Atiyah-Singer density for $g(\ep)$ converges to that for $g_\Phi$, i.e.\ 
\begin{equation}\label{eq:limit-of-as}
\lim_{\ep\searrow 0}\int_X \mbox{AS} (g(\ep))\;=\;\int_X \mbox{AS}(g_\Phi)\,.
\end{equation}
\end{proposition}
\begin{proof}
This is a computation. We shall use the method of moving frames, cf.\  
\cite{Spivak} for more on this formalism. Recall that if $\omega^0, \ldots, \omega^n$
is any orthonormal set of one-forms, then the connection one-forms $\omega_i^{\, j}$ 
are determined uniquely by the equations
\[
d\omega^i = \omega^j \wedge \omega_j^{\, i}, \qquad \omega_j^{\, i} = -\omega_i^{\, j}.
\]
From these we define the curvature two-forms
\[
\Omega_i^{\, j} = d\omega_i^{\, j} - \omega_i^{\, k}\wedge \omega_k^{\, j}.
\]
Here, and elsewhere below, summation on repeated indices is intended.

The strength of this method, of course, is that it can be adapted to the specific
geometry. Thus here we shall choose the coframe for $g(\ep)$ as follows. 
Let $Y=\pa X$. Choose an orthonormal coframe $\tilde{\omega}^\alpha$, $1 \leq \alpha \leq k$, for $(Y,h)$, 
and $\omega^\mu$, $k+1 \leq \mu \leq n$, for the restriction of $\kappa$ to each fibre. 
These forms may also depend smoothly on $\ep$ and $x$ (in $x\geq 0$, $\ep \geq 0$),
and in addition, the $\omega^\mu$ may also depend on $y \in Y$.
In the following, we shall use the Chern convention that Roman indices $i,j, \ldots$ vary
between $0$ and $n$, while the Greek indices $\alpha, \beta, \ldots$ vary between
$1$ and $k$ and $\mu, \nu, \ldots$ vary between $k+1$ and $n$. Now define
\[
\omega^0 = \frac{dx}{x(x+\ep)}, \qquad \omega^\alpha = \frac{\phi^*(\tilde{\omega}^\alpha)}{x+\ep};
\]

Then 
\[
\{\omega^0, \omega^1, \ldots, \omega^k, \omega^{k+1}, \ldots, \omega^n\}
\]
is an orthonormal coframe for $g(\ep)$. 

After some computation we obtain
\begin{eqnarray*}
d\omega^0  &=& 0 \\
&\equiv & \omega^\alpha \wedge \omega_\alpha^{\, 0} + \omega^\mu \wedge \omega_\mu^{\, 0}
\\d\omega^\alpha &=& -\frac{dx}{(x+\ep)^2} \wedge \phi^*(\tilde{\omega}^\alpha) + 
\frac{dx}{(x+\ep)} \wedge (\phi^*(\tilde{\omega}^\alpha))^{\prime} +
\omega^\beta \wedge \phi^*(\tilde{\omega}_{\beta}^{\, \alpha}) \\
&\equiv& \omega^0 \wedge \omega_0^{\, \alpha} + \omega^\beta\wedge \omega_{\beta}^{\, \alpha} + \omega^\mu
\wedge \omega_{\mu}^{\, \alpha} \\
d\omega^\mu &=& dx \wedge (\omega^\mu)^\prime + (x+\ep) \omega^\alpha\wedge E_{\alpha}^{\, \mu} +
\omega^\nu \wedge E_\nu^{\, \mu} \\
&\equiv& \omega^0 \wedge \omega_0^{\, \mu} + \omega^\alpha \wedge \omega_\alpha^{\, \mu} +
\omega^\nu \wedge \omega_\nu^{\, \mu}.\\
\end{eqnarray*}
Here the ${}^{\prime}$ denotes differentiation with respect to $x$, and $E_i^{\, j}$ denotes terms 
(involving the curvature and second fundamental form of the fibres) which are uniformly bounded
(with respect to the unscaled metric $\tilde{g}$ on $\pa X$) along with their derivatives as 
$x,\ep \to 0$. 

More specifically, in the formula for $d\omega^\alpha$, we use that $d$ commutes with $\phi^*$. 
The expression for $d\omega^\mu$ contains $\omega^\alpha$ factors corresponding
to the derivative of the fibre metric in the horizontal direction, and also to the
variation of the horizontal subspaces in the fibre direction. We refer to
\cite{HHM} \S 5.3.1 (particularly (43)-(45)) for the precise details, but note
simply that the $\omega^\alpha \wedge \omega^\nu$ components correspond to the 
second fundamental form in the normal direction $e_\alpha$ to the fibre (with 
respect to the scaled metric on $\pa X$ for a given $x$ and $\ep$), and are indeed
of the form $(x+\ep)E_\alpha^{\, \mu}$, while the $\omega^\alpha \wedge \omega^\beta$
components in $d\omega^\mu$ correspond to the curvature of the horizontal distribution, which are 
of the form $(x+\ep)^2 E_\alpha^{\, \mu}$, hence even lower order. Next, the terms $E_\nu^{\, \mu}$
are precisely the connection one-forms $\omega_\nu^{\, \mu}$ for the metric induced by $\kappa$ on
the fibres; in particular, these do not involve any $\omega^\alpha$ factors.
Finally, we have included the extra terms involving the $x$ derivative of $\tilde{\omega}^\alpha$ and
$\omega^\mu$ since we do allow the metric $h$ on $Y$ and symmetric two-tensor $\kappa$ 
to depend smoothly on $x$. 

Using this same $E_i^{\, j}$ notation for all `negligible' bounded terms, we now claim that
\[
\omega_0^{\, \alpha} =  -\frac{x}{x+\ep} \tilde{\omega}^\alpha +  x E_0^{\, \alpha},  \qquad 
\omega_0^{\, \mu} =  x(x+\ep) E_0^{\, \mu}, 
\]
\[
\omega_\alpha^{\, \beta} = E_\alpha^{\, \beta},  \qquad  \omega_\alpha^{\, \mu} = (x+\ep)E_\alpha^{\, \mu},
\qquad \omega_\mu^{\, \nu} = E_\mu^{\, \nu}.
\]

To verify this, we simply need to show that these forms satisfy the structure equations
and are skew-symmetric in their indices, for then Cartan's lemma guarantees uniqueness.
The equations for all terms except the $\omega_\alpha^{\, \mu}$ (which by skew-symmetry,
we {\it require} to be equal to $-\omega_\mu^{\, \alpha}$) are clear enough.
For these terms, first note that the equation for $d\omega^\alpha$ has no vertical
components, which means that $\omega^\mu \wedge \omega_\mu^{\, \alpha}$ must vanish.
This means that 
\[
\omega_\mu^{\, \alpha} = c_{\alpha,\mu,\nu} \omega^\nu, \qquad \mbox{and}\qquad
c_{\alpha,\mu,\nu} = c_{\alpha,\nu,\mu}.
\]
(The point is that there can be no $\omega^\beta$ or $\omega^0$ components.)
Finally, setting this into the equation for $d\omega^\mu$, and noting that the 
$E_\nu^{\,\nu}$ term is already accounted for by the $\omega_\nu^{\,\mu}$, we
must have $\omega_\alpha^{\,\mu} = (x+\ep)E_\alpha^{\,\mu}$, as claimed.

When computing each of the curvature two-forms $\Omega_i^{\, j}$, we write all forms in terms of 
$dx$, $\tilde{\omega}^\alpha$ and $\omega^\mu$, which are smooth in the ordinary sense
up to $\ep = x=0$. We single out the particular terms which help or hurt us, and as before 
gather all the harmless remaining factors into terms $F_i^{\, j}$, which are 
uniformly bounded in $x,\ep \geq 0$. Thus, after further work, we obtain
\begin{eqnarray*}
\Omega_0^{\, \alpha} & = & dx \wedge \left(-\frac{\ep}{(x+\ep)^2} \tilde{\omega}^\alpha
+ F_0^{\, \alpha}\right) + x F_0^{\, \alpha}, \\
\Omega_0^{\, \mu} & = & (x+\ep) F_0^{\, \mu}, \\
\Omega_\alpha^{\, \beta}  &=& F_{\alpha}^{\, \beta}, \\
\Omega_\alpha^{\, \mu} &=& dx \wedge F_\alpha^{\, \mu} + (x+\ep) F_\alpha^{\, \mu}, \\
\Omega_\mu^{\, \nu} &=& F_{\mu}^{\, \nu}.
\end{eqnarray*}
Only the first of these requires more explanation. We have
\[
d\omega_0^{\, \alpha} - \omega_0^{\, \beta}\wedge \omega_\beta^{\, \alpha} - \omega_0^{\, \mu}\wedge \omega_\mu^{\,\alpha}
\]
\[
= d\left( -\frac{x}{x+\ep}\right)\wedge \tilde{\omega}^\alpha + dx \wedge F_0^{\, \alpha} - \frac{x}{x+\ep}
\left(d\tilde{\omega}^\alpha - \tilde{\omega}^\beta \wedge \tilde{\omega}_\beta^{\, \alpha} \right) + x F_0^{\, \alpha}.
\]
The first terms on the right, involving $dx$, and the final term, correspond to the assertion above. The
middle terms appear not to be of the correct form, but the particular combination
in parentheses is just the structure equation for the connection one-forms and hence vanishes.

Recalling that $\dim X = n+1 = 2 m$, the integral $\int_X \widehat{A}(g(\epsilon))$ is a 
linear combination of terms of the form:
\[
\int_X {\rm Tr}\, R^{m_1}(\epsilon ) \ldots {\rm Tr}\, R^{m_p}(\epsilon), \qquad  m_1 +\ldots +m_p=m.
\]
To fix the ideas and simplify the notation we focus on 
\[
\int_X {\rm Tr}\, R^{m}(\epsilon ),
\]
since all other terms are handled the same way. In terms of the curvature two-forms,
\begin{equation} \label{eq:tr}
{\rm Tr} R^m(\epsilon)\,=\, \sum 
\Omega^{\, i_2}_{i_1} \Omega^{\, i_3}_{i_2} \cdots \Omega^{\, \, i_{1}}_{i_m}.
\end{equation}

Now substitute in this the expressions we have obtained for the $\Omega_i^{\, j}$. 
Using the boundedness of all of the $E_{i}^{\, j}$, the only terms in any of 
these curvature forms which is not bounded near $x=\ep=0$ is $\Omega_0^{\,\alpha}$, 
and in fact only its first term $\ep (x+\ep)^{-2}dx \wedge \tilde{\omega}^\alpha$ causes difficulties.
Thus we may as well suppose that this is the first term, i.e.\ $i_1 = 0$ and 
$i_2 = \alpha$, and we can replace the entire two-form $\Omega_{0}^{\, \alpha}$
by this single bad term. The final term in the entire product is either $\Omega_\mu^{\, 0}$ or 
$\Omega_\beta^{\, 0}$ for some $\mu$ or $\beta$. In the former case this contains a vanishing 
factor $(x+\ep)$, while in the latter, only the part of this two-form which  
does not contain a $dx$ contributes, and this has the same vanishing factor.
Thus in all cases, the entire $2m$-form is bounded (though not necessarily smooth!)
near $\ep = x = 0$, and we can pass to the limit, as desired.
\end{proof}

\subsection{Adiabatic limit of the eta invariant}\label{subsec:eta}

We briefly recall the context of the Bismut-Cheeger theorem \cite{BC}. 
Let $M$ be an odd dimensional, compact spin manifold which is 
the total space of a fibration 
\[
Z\rightarrow M\xrightarrow{\phi} Y
\]
where the base $Y$ is also spin. We fix a connection $TM=T_H (M)\oplus T(M/Y)$,
where $T_H (M)\simeq \phi^* (TY)$ and $T(M/Y)$ denotes the vertical tangent bundle. 
Let $h$ be a Riemannian metric on $Y$ and $\kappa$ a symmetric two-tensor on $TM$ 
which restricts to a metric on each $Z_y$ and which annihilates the horizontal space,
and introduce the Riemannian submersion metric $\tilde{g} :=\phi^* h + \kappa$.
Denote by $\nabla$ and $\nabla^{M/Y}$ the Levi-Civita connection for $g_M$ 
and the induced connection on $T(M/Y)$ obtained by compressing $\nabla$ by 
the projections $P:TM\to T(M/Y)$. Let $S$ be the vertical spinor bundle and 
$E\to M$ an additional Hermitian bundle endowed with a unitary connection.
The bundle $F:=S\otimes E$ is a vertical Clifford module. 
Finally, let $\FF:= \phi^* (\Lambda^* Y)\otimes F$. To fix the notation we assume 
that the fibers are even dimensional. 

To this entire set of data one associates the rescaled Bismut superconnection
\[
\AAA_t: C^\infty (M,\FF)\rightarrow  C^\infty (M,\FF)\,,
\]
cf.\ \cite{Bismut-inventiones} and \cite{BGV}. The operator $d\AAA_t/dt \exp (-\AAA^2_t)$ 
is a vertical family of smoothing operators $(\mathcal{K}_y)_{y\in Y}$
with coefficients which are differential forms on the base $Y$. From this 
family one obtains a differential form of odd degree on the base $Y$, 
$${\rm Str}\left( \frac{d\AAA_t}{dt} \exp (-\AAA^2_t)\right).$$
The value of this form at $y\in Y$ is obtained by
restricting $\mathcal{K}_y$ to  the diagonal 
$\Delta_y\subset \phi^{-1}(y)\times \phi^{-1}(y)$, taking its supertrace and then 
integrating over $\Delta_y$. 

Assume now that the vertical family of Dirac operators $(D_y)_{y\in Y}$ associated to 
the data above is invertible. Then the integral
$$
\int_0^\infty \frac{1}{\sqrt{\pi}} {\rm Str}\left( \frac{d\AAA_t}{dt} \exp
(-\AAA^2_t)\right) dt
$$
converges and defines the eta form $\widetilde{\eta}\in C^\infty (Y,\Lambda^*Y)$
associated to the family $(D_y)_{y\in Y}$.

The adiabatic limit formula of Bismut and Cheeger states that if 
$\eta(\ep)$ is the eta invariant for the Dirac operator associated
to the metric 
$$
g_M(\ep):= \frac{\phi^* h}{\epsilon^2}+\kappa\,;$$
then
$$
\lim_{\ep\to 0} \eta(\ep)=\int_Y \widehat{A}(Y,h)\wedge
\widetilde{\eta}\,.
$$

Applied to the boundary operator $D^\pa_{g(\ep)}$ on $M = \pa X$, 
we obtain the limiting behaviour of the final term in (\ref{eq:phigb2}).
This completes the proof of the index formula for Dirac operators
associated to $\Phi-$metrics in the fully elliptic case.
{\small
\bibliographystyle{plain}
\bibliography{phi-index}

\begin{thebibliography}{10}

\bibitem{BGV}
Nicole Berline, Ezra Getzler, and Mich{\`e}le Vergne.
\newblock {\em Heat kernels and {D}irac operators}, volume 298 of {\em
  Grundlehren der Mathematischen Wissenschaften [Fundamental Principles of
  Mathematical Sciences]}.
\newblock Springer-Verlag, Berlin, 1992.

\bibitem{Bismut-inventiones}
Jean-Michel Bismut.
\newblock The {A}tiyah-{S}inger index theorem for families of {D}irac
  operators: two heat equation proofs.
\newblock {\em Invent. Math.}, {\bf 83}(1):91--151, 1985.

\bibitem{BC}
Jean-Michel Bismut and Jeff Cheeger.
\newblock {$\eta$}-invariants and their adiabatic limits.
\newblock {\em J. Amer. Math. Soc.}, {\bf 2}(1):33--70, 1989.

\bibitem{HHM}
Tam{\'a}s Hausel, Eugenie Hunsicker, and Rafe Mazzeo.
\newblock Hodge cohomology of gravitational instantons.
\newblock {\em Duke Math. J.}, {\bf 122}(3):485--548, 2004.

\bibitem{LM0}
Robert Lauter and Sergiu Moroianu.
\newblock Homology of pseudodifferential operators on manifolds with fibered
  cusps.
\newblock {\em Trans. Amer. Math. Soc.}, {\bf 355}(8):3009--3046 (electronic),
  2003.

\bibitem{LM}
Robert Lauter and Sergiu Moroianu.
\newblock An index formula on manifolds with fibered cusp ends.
\newblock {\em J. Geom. Anal.}, {\bf 15}(2):261--283, 2005.

\bibitem{Ma-Mel-phi}
Rafe Mazzeo and Richard~B. Melrose.
\newblock Pseudodifferential operators on manifolds with fibred boundaries.
\newblock {\em Asian J. Math.}, {\bf 2}(4):833--866, 1998.
\newblock Mikio Sato: a great Japanese mathematician of the twentieth century.

\bibitem{mel-rochon}
Richard Melrose and Fr\'ed\'eric Rochon.
\newblock Index in k-theory for families of fibred cusp operators.
\newblock Preprint 2005.

\bibitem{Melrose}
Richard~B. Melrose.
\newblock {\em The {A}tiyah-{P}atodi-{S}inger index theorem}, volume~4 of {\em
  Research Notes in Mathematics}.
\newblock A K Peters Ltd., Wellesley, MA, 1993.

\bibitem{phi-vs-di}
Sergiu Moroianu.
\newblock Fibered cusps versus $d$-index theory.
\newblock preprint 2006.

\bibitem{Mueller}
Werner M{\"u}ller.
\newblock {\em Manifolds with cusps of rank one}, volume~{\bf 1244} of {\em
  Lecture Notes in Mathematics}.
\newblock Springer-Verlag, Berlin, 1987.
\newblock Spectral theory and $L\sp 2$-index theorem.

\bibitem{Nye-Singer}
Tom M.~W. Nye and Michael~A. Singer.
\newblock An {$L\sp 2$}-index theorem for {D}irac operators on {$S\sp
  1\times\bold R\sp 3$}.
\newblock {\em J. Funct. Anal.}, {\bf 177}(1):203--218, 2000.

\bibitem{Spivak}
Michael Spivak.
\newblock {\em A comprehensive introduction to differential geometry. {V}ol.
  {II}}.
\newblock Publish or Perish Inc., Wilmington, Del., second edition, 1979.

\bibitem{Vaillant-phd}
Boris Vaillant.
\newblock Index- and spectral theory for manifolds with generalized fibred
  cusps, 2001.
\newblock Dissertation, Rheinische Friedrich-Wilhelms-Universit\"at Bonn, Bonn,
  2001.

\end{thebibliography}

\end{document}